\documentclass[a4paper, 10pt, conference]{ieeeconf}
\IEEEoverridecommandlockouts 
\usepackage{graphicx}
\usepackage[intlimits]{amsmath}
\usepackage{amsfonts}
\usepackage{amssymb}
\usepackage{hhline}

\newcommand{\be}{\begin{equation}}
\newcommand{\ee}{\end{equation}}

\newcommand{\trn}{^{\rm\scriptscriptstyle T}}

\newcommand{\tet}{\vartheta}

\newcommand{\eps}{\varepsilon}
\newcommand{\ph}{\varphi}
\newcommand{\mR}{\mathbb R}
\newcommand{\inr}{\!\in \mR}

\newtheorem{theorem}{Theorem}

\def\CheckPDFoutput{%
\CheckPDFoutput%
\ifx\unprotect\undefined%
 \DeclareGraphicsRule{.jpg}{bmp}{}{}%
  \else%
  \pdfoutput=1%
\fi%
}
\title{Adaptive Observer-based Synchronization of Nonlinear Nonpassifiable
Systems}
\author{V.O.~Nikiforov,\thanks{{T}{his} work was partly supported by Russian Foundation
of Basic Research (RFBR), grant 05-01-00869 and Scientific Program of RAS No 19.
}\thanks{V.O.~Nikiforov is with the
St.Petersburg State University of Fine Mechanics and
 Optics, Saint Petersburg, Russia. E-mail:
 nikiforov@mail.ifmo.ru.}  A.L.~Fradkov, and B.R.~Andrievsky\thanks{A.L.~Fradkov and B.R.~Andrievsky
are with the Control of Complex Systems Lab., Institute for
Problems of Mechanical Engineering of
 Russian Academy of Sciences (IPME RAS), St.Petersburg,
Russia. E-mail: \{alf,bandri\}@control.ipme.ru;
 bandri@yandex.ru .}
}

\begin{document}
%\markboth{IEEE Transactions On Automatic Control, Vol. XX, No. Y,
%Month 2003}
%{Murray and Balemi: Using the style file IEEEtran.sty} %!PN
%{Nikiforov {\it et al.}: Adaptive Observer-based Synchronization} %!PN

\maketitle

\begin{abstract}
In this paper the relative degree limitation for adaptive
observer-based synchronization schemes is overcome. 
The scheme is extended to nonpassifiable systems. Two 
synchronization methods are described and 
justified based on augmented error adaptive 
observer and high-order tuners. 
The solution is based on modern theory of 
nonlinear adaptive control, particularly on 
nonlinear observer structure and new classes of adaptation algorithms.
Conditions of parametric convergence of
the parameter estimation are established for the noiseless case.
Robustness of the scheme to the bounded measurement error is
established. The results are illustrated by example of application the
proposed adaptive synchronization of chaotic Lorenz systems.
\end{abstract}


\section{Introduction}\label{Sec:1}

\PARstart{S}{ynchronization} has fo\-und va\-ri\-ous ap\-pli\-ca\-tions 
during last decade. Particularly, a lot of
interest has been attracted to the problem of information transmission by means 
of chaotic signals modulation, see 
\cite{IEEE97,IEEE01,CTA99,KK99}. A number of results in 
this area are based on adaptive synchronization approach
\cite{F95,MF96,YC96,FM97}. New method of adaptive synchronization
exploiting Lyapunov functions and passification was developed in
\cite{F95,FM97}
and later extended to observer-based adaptive synchronization withapplication to telecommunication \cite{FNM00,AF00}. The possibility of fast 
transmitting messages in noisy channel using new approach
has been demonstrated \cite{A00,A02}. However, applicability of the
method, proposed in \cite{F95,FM97} is restricted by the passifiability of the
plant (the master system).
It implies the relative degree limitation: the relative
degree of the plant model should be equal to zero or one. Such a
limitation prevents from increasing security of communications by using
master system with higher relative degree \cite{Hui00}.

In this paper we overcome the relative degree limitation for adaptive
observer-based synchronization schemes \cite{NM97} and extend them to
nonpassifiable
systems, particularly to systems with relative degree greater than one.
The solution is based on modern theory of nonlinear adaptive control
\cite{FMN99,KKK,M90}, particularly on nonlinear observer structure and
new classes of adaptation algorithms \cite{NikAut,NV01}.

In the literature the methods of adaptive synchronization
suitable for adaptive systems with relative degree greater than one were 
proposed \cite{Hui00,Ge00}. However, 
the algorithm of \cite{Ge00} provides state 
feedback rather than output feedback and does 
not allow for presence of unknown 
parameters in the master system. The
approach of \cite{Hui00} is based on canonical forms for linear
adaptive
observers (see, e.g.  \cite{SB89}).
Design of adaptive observer for nonlinear systems in \cite{Hui00}
requires
case by case consideration and special tricks for each nonlinear
system.
For example, it is not clear how to apply results of \cite{Hui00} to
Lorenz
systems. Finally, the results of \cite{Hui00} do not allow to
incorporate
measurement errors.

Unlike \cite{Hui00,Ge00}, in this paper a unified approach based on
scheme of  \cite{AF00,F95,FM97,FNM00} and new classes of adaptation
algorithms
of \cite{FMN99,NikAut,NV01} are proposed. It allows to cope with
bounded noise
by means of robust modification of adaptation algorithm.
Simulation results demonstrate better convergence and
robustness properties the adaptive observer (slave system)
compared with results of \cite{AF00}.

In Sections \ref{Sec:2}, \ref{Sec:3} general method of adaptive
observer design
for higher relative degree systems is described and justified. To
clarify the
presentation of main idea we start with the simple disturbance-free
case and
introduce basic schemes of adaptive observers in Section \ref{Sec:2}.
In  Sec. \ref{Sec:3} we modify basic schemes to provide
robustness in the presence of external disturbances (measurement
noise). The method is applied to Lorenz system in
Section \ref{Sec:4} where simulation results conforming the theoretical
statements are presented. Numerical examples demonstrating applicaton of 
the proposed scheme to signal transmission are given. Preliminary version 
of the results was announced in \cite{CDC02}.

\section{Problem Statement}

{F}{ollowing} \cite{FNM00} we assume that the plant (the {\it
master system}) is
described by the state-space equations of the form:
\begin{align}
\dot x&=Ax+\ph_0(y)+b\ph\trn(y)\theta,~~~ y=c\trn x,
\label{2}\\
y_r&=y+\xi,\label{27}
\end{align}
where $x\in\mR^n$ is the inaccessible state vector, $y$ is
the plant output (transmitted signal), $y_r$ is the measurable noisy
signal,
$\xi$ is the additive channel noise (presented by a bounded function of
time), $\theta\in\mR^m$ is the unknown vector of the plant parameters
(possibly representing a message). It is assumed that the
nonlinearities
$\ph_0(y)$, $\ph(y)$, matrix $A$ and vectors $b$, $c$ are known.

 We accept the following plant-model assumptions.

{\bf Assumption 1:} { for any bounded initial condition $x(0)$ and any
value of
vector $\theta$ the
state $x(t)$ is a bounded function of time.}

{\bf Assumption 2:} {functions $\ph_0(y)$, $\ph(y)$ are bounded for any
bounded
$y$.}

The problem is to design an adaptive observer (a dynamical system) of
the form
\begin{equation}
\dot z=F(z,y_r),\quad\hat\theta=h(z,y_r),
\label{2z}
\end{equation}
such that
$\overline{\lim\limits_{t\to\infty}}|\theta-\hat\theta|\le\Delta$,
where $\Delta$ is some positive constant.

\section{Design of Adaptive Observers: Disturban\-ce\--Free Case}\label{Sec:2}

In this section we assume that $\xi(t)\equiv 0$. Then the
adaptive observer (the {\it slave system}) may have the
following structure:
\begin{align}
\dot{\hat x}&=A\hat x+\ph_0(y)+b\ph(y)\trn\hat\theta+
k(y-\hat y),\quad 
\hat y=c\trn\hat x, \label{eq:3}\\
\dot{\hat\theta}&=F_\theta\big(\hat\theta,\hat x,y\big),
\label{2o}
\end{align}
where $\hat x$ and $\hat y$ are the estimates of $x$ and $y$, $\hat
\theta$
is the vector of adjustable parameters (representing estimate of the
parameter
$\theta$), and $z={\rm col}(x,\hat\theta)$.

Such a structure was proposed in \cite{FNM00} where convergence
conditions were established for the case of plant (master system) with 
passifiable linear part (for bounded
disturbances see \cite{AF00}). Passifiability condition imposes strong
restriction on the plant model: its relative degree $r$ should be equal
to $1$.
At the same time it would be important to have extended solutions to
this
problem in the case $r>1$.

Below we present two solutions to the posed problem under assumptions
accepted.
The first
solution is based on the use of the {\it augmented error} (AE) concept
\cite{FMN99,N&F94,N&A89}, while the second one utilizes the idea of the
{\it
high-order tuners}
(HOT)
\cite{FMN99,NikAut,NV01,Mor92}.

To employ the augmented error concept we use the adaptive observer
\eqref{eq:3},
\eqref{2o} where vector $k$ is chosen so that $F=A-kc\trn$ is Hurwitz.
To derive the adaptation algorithm, we first obtain so-called
{\it error model}. Differentiating estimation error $\eps=x-\hat x$ in
view of (\ref{2}) and \eqref{eq:3}, we obtain
\begin{equation}
\dot\eps =F\eps+b\ph(t)\trn\tilde\theta, \quad e=c\trn\eps \label{5}
\end{equation}
where $\ph(t)=\ph\big(y(t)\big)$,
$\tilde \theta=\theta-\hat\theta$ is the parameter error, while
$e=y-\hat y$ is the output error accessible to measurements. The error
model (\ref{5}) can be rewritten in the form:
\begin{equation}
e=H(p)\Big[\ph(t)\trn\tilde \theta\Big],
\end{equation}
where $p=d/dt$ is the differential operator, and
the transfer function $H(p)=c\trn(p{\rm I}-F)^{-1}b$ is asymptotically
stable.

The adaptation algorithm can be chosen in the form
\cite{FMN99,N&A89,N&F94}:
\begin{equation}
\dot{\tilde  \theta}=\gamma\omega(t)\trn\hat e,
\label{7}
\end{equation}
where $\gamma>0$ is the adaptation gain, $\omega(t)=H(p)[\ph(t)]$ is
the filtered regressor, while the {\it augmented error} $\hat e$ is
defined by
the following equation:
\begin{equation}
\hat e=e+H(p)\Big[\ph(t)\trn\hat \theta\Big]-\omega(t)\trn\hat \theta.
\label{8}
\end{equation}

%%%%%%%%%%%%%%%%%%%%%%%%%%%%%%%%
Introduce the following definition. 

\par{\it Definition:}
A vector-function $f:[0,\infty) \to \mR^{m}$
is called {\it persistently exciting (PE) on } $[0,\infty)$, if it is
measurable and bounded on $[0,\infty)$ and there exist $\alpha >0, T>0$
such that
\begin{equation}
\int_{t}^{t+T }f(s)f(s)^{T}ds \ge\alpha I
\label{c8}
\end{equation}
for all $t\ge 0$.

%%%%%%%%%%%%%%%%%%%%%%%%%%%%%%

\begin{theorem}\label{th:1}
The closed-loop system consisting of the master system
(\ref{2}), adjustable observer \eqref{eq:3} and algorithm
of adaptation (\ref{7}), (\ref{8}) has the following properties:
\begin{enumerate}
\item[{\it i})] for any initial conditions and any $\gamma>0$ all the
closed-
loop
signals are bounded and
\begin{equation}
\lim\limits_{t\to\infty} \left(y(t)-\hat y(t)\right)=0;
\label{u1}
\end{equation}

\item[{\it ii})] if the vector function $\ph(t)$ satisfies PE condition
and the
transfer function $H(p)$ is minimum phase then in
addition to ({\it i})
\begin{equation}
\lim\limits_{t\to\infty} \big|\theta-\hat \theta(t)\big|=0.
\label{u2}
\end{equation}
\end{enumerate}
\end{theorem}
\begin{proof}
It is known that for the augmented error $\hat e$ \eqref{8} one can
write the following equivalent model \cite{FMN99,N&A89}, neglecting 
exponentially vanishing term due to nonzero initial conditions:
\begin{equation}
\label{pr:A}
\hat e =\omega(t)\trn\tilde\theta.
\end{equation}
Then differentiating Lyapunov function
$V(\tilde\theta)=\dfrac{1}{2\gamma}\tilde
\theta\trn\tilde\theta$ along solutions of \eqref{7} and in view of
\eqref{pr:A}, we obtain $\dot V(\tilde\theta)=-\hat{e}^2$. The latter
means
boundedness of $\hat\theta$ and zeroing $\hat e(t)$ (since the
right-hand sides
of  \eqref{2},\eqref{eq:3} and \eqref{7} are locally Lipshitz in $x$,
$\hat x$ and
$\hat \theta$ uniformly in $t$ \cite{KKK}). Since $\omega(t)$ is
bounded, from
\eqref{7} we have that $\dot{\hat\theta}\to 0$ as $t\to \infty$.
Therefore, from \eqref{8} we obtain that $e -\hat e\to 0$ and, as a
consiquence,
$e\to 0$. Part ({\it i}) is proved. Part ({\it ii}) can be
straightforwardly
proved with the use of standard arguments \cite{N&A89}.
\end{proof}

Now we present an alternative solution to the posed problem utilizing
an idea of the high-order tuners. For this, we use the following
adjustable observer:
\begin{equation}
\dot{\hat x}=A\hat x+\ph_0(y)+b\nu(y,\hat \theta)+k(y-\hat y)
,~~~~\hat y= c\trn\hat x \label{10}
\end{equation}
where the adjustable feedback $\nu(y,\hat \theta)$ will be defined
below.
In this case the error model takes the view
$$ \dot\eps =F\eps +b(\ph (t) \trn \theta-\nu),\quad e=c\trn\eps $$
or
\begin{equation}
 e=H(p)\Big[\ph(t)\trn\theta-\nu\Big].\label{10a}
\end{equation}
Chose a transfer function $W(p)$ obeying the equation:
$$ W(p)=(p+\lambda)H(p)$$
where $\lambda$ is any positive constant. Then the model (\ref{10a})
can be
rewritten in the form
\begin{equation}
e=\dfrac{1}{p+\lambda}\Big[\varpi(t)\trn\theta-W(p)[\nu]\Big]
\label {11}
\end{equation}
where $\varpi(t)=W(p)[\ph(t)]$. Analysis of the model (\ref{11})
motivates the following choice of the adjustable feedback:
\begin{equation}
\nu=W(p)^{-1}\Big[\varpi(t)\trn\hat \theta\Big]
\label{12}
\end{equation}
where $\hat \theta$ is the vector of adjustable parameters. To realize
the
feedback (\ref{12}) we need to generate not only the adjustable
parameters
$\hat \theta$, but also their high-order derivatives up to
order $r-1$, where $r$ is the relative degree of the transfer
function $H(p)$. To overcome this problem we use the following
algorithm
of adaptation:
\begin{align}
\dot{\hat \psi}_i&=\varpi_i e,
\label{13}
\\
\dot\eta_i&=(1+\mu\varpi\trn\varpi)(\Gamma\eta_i+h\hat \psi_i),
\label{14}\\
\hat \theta_i&=l\trn\eta_i
\label{15}
\end{align}
where $i=1, 2, \ldots, m$, $\mu>0$ is a design parameter, and
$(l,\Gamma,h)$ is a minimal realization of the transfer function
$\alpha(0)/\alpha(p)$ with a Hurwitz polynomial $\alpha(p)$ of order
$r-2$,
i.e. $\alpha(0)/\alpha(p)=l\trn(p{\rm I}-\Gamma)^{-1}h$.

{\it Remark.} If $r\le 2$, one needs not to use 
additional filters \eqref{14},\eqref{15}. In this particular case the
algorithm of adaptation takes the form $\dot{\hat {\theta}}=\varpi e$.

\begin{theorem}\label{th:2}
The closed-loop system consisting of the master system
(\ref{2}), observer (\ref{10}) and
adjustable feedback (\ref{12})--(\ref{15}) has the following
properties:

\par{\it i}) for any initial conditions and any
$$\mu>\dfrac{3}{4\lambda}\Big(|l|+|PF^{-1}h|\Big)^2$$
where the positive definite matrix $P$ obeys the equality
$F\trn P+PF=-2{\rm I}$, all the closed-loop signals are bounded and
regulation (\ref{u1}) is achieved;

\par{\it ii}) if the vector function $\ph(t)$ satisfies PE condition
and
the transfer function $H(p)$ is minimum phase then
in addition to ({\it i}) asymptotic convergence (\ref{u2})
is guaranteed.
\end{theorem}

{\it Proof} of the theorem is based on the standard arguments which can
be
found, for example, in \cite{FMN99,NikAut}.

Finally, consider a more general case when the master system is
described by the
following equations:
\begin{equation}\dot x= A(y)x+\ph_0(y)+b\ph(y)\trn\theta,\quad
y=c\trn x.
\label{20}\end{equation}
We  accept the following additional plant model assumptions.

{\bf Assumption 3:} There exist a vector function $k(y)\in\mR^n$ and
scalar function $V(x)$ such that
$$
\begin{array}{l}
c_1|x|^2\le V(x) \le c_2|x|^2, \vspace{2mm} \\
\dfrac{\partial V}{\partial x}(x)\left(A(c\trn x)-k(c\trn
x)c\trn\right)x
\le -c_3|x|^2, \vspace{2mm} \\
\left|\dfrac{\partial V}{\partial x}\right|\le c_4|x|,
\end{array}
$$
where $c_i$ are some positive constants ($i=1$, $2$, $3$, $4$).

{\bf Assumption 4:} All entries of matrix $A(y)$ are bounded for any
bounded
$y$.

In other words, Assumption 3 means that the autonomous system
$$\dot x=G(c\trn x)x,$$
where $G(c\trn x)=G(y)=A(y)-k(y)c\trn$, is exponentially stable
\cite{Kras63}.

It is worth noting that the design methods presented above are not
applicable to the model (\ref{20}). Therefore we introduce
one more method presenting a special kind of the scheme with augmented
error.
For this we employ the following adjustable observer
\begin{eqnarray}
\dot{\hat x}=&A(y)\hat x+\ph_0(y) +b\ph(y)\trn\hat \theta
+k(y)(y-\hat y), \label{22} \\
& \qquad\hat y=c\trn\hat x\nonumber 
\end{eqnarray}
where the time-varying vector $k(y)$ is chosen so that Assumption 3 is
valid. In this case the error model takes the view:
\begin{equation} \dot \eps =G(t)\eps +b\ph(t)\trn\tilde \theta,\quad
e=c\trn \eps,
\label{nA}\end{equation}
where $G(t)=G(c\trn x(t))$.

Define the augmented error as follows:
\begin{equation}
\hat e=e+c\trn\eta,
\label{23}
\end{equation}
where the auxiliary vector $\eta$ is generated by the filters:
\begin{eqnarray}
\label{24}
\begin{aligned}
\dot \eta\!=\!G(t)\eta - \Omega\dot{\hat \theta},
\quad\eta\in\mR^n,
\\
\dot \Omega \!=\! G(t)\Omega+b\ph(t)\trn,\quad\Omega\in\mR^{n\times n}.
\label{25}
\end{aligned}
\end{eqnarray}
Then adaptation algorithm can be chosen in the form \cite{N&F94}:
\begin{equation}
\dot{\hat \theta}=\gamma \omega\trn\hat e,
\label{26}
\end{equation}
where $\omega=c\trn\Omega$.

\begin{theorem}\label{th:3}
The closed-loop system consisting of the master system
(\ref{20}), adjustable observer (\ref{22}),
scheme of augmentation (\ref{23})--(\ref{25}) and algorithm
of adaptation (\ref{26}) has the following properties:
\par{\it i}) for any initial conditions and any $\gamma>0$ all the
closed-loop signals are bounded and regulation (\ref{u1}) is
achieved;

\par{\it ii})  if, in addition, the vector function $\ph(t)$ satisfies
PE
condition then asymptotic convergence (\ref{u2}) is guaranteed.

\end{theorem}
\begin{proof}
Differentiating the following auxiliary error
$\delta=\eps+\eta-\Omega\tilde\theta$ in view of
equations \eqref{nA} and \eqref{24} we obtain
\begin{align*}
\dot\delta&=G\eps+b\ph\trn\tilde\theta+G\eta-\Omega\dot{\hat\theta}
-G\Omega\hat\theta-b\ph\trn\tilde\theta+\Omega\dot{\hat\theta}\\
&= G\big(\eps+\eta-\Omega\tilde\theta\big)=G\delta.
\end{align*}
Then for the augmented error defined by equation \eqref{23} we can
write the following equivalent model: $\hat e=\omega\trn\tilde\theta+\delta_e$,
where $\delta_e=c\trn \delta$ exponentially vanishes. Then using the
same
arguments as in the proof of Theorem \ref{th:1} we can show boundedness
of all the closed-loop signals, regulation \eqref{u1} and convergence 
\eqref{u2} (under PE condition).

\end{proof}

\section{Robust Adaptive Obser\-vers: Noisy
Measurements}\label{Sec:3}

{I}{t} is known that adaptation algorithms of pure integral
action can
loss stability in the presence of external disturbances or noise of
measurements   \cite{FMN99,I&K84}. In this section we present
robustified
modifications of above schemes which are applicable in the case of
noisy
measurements.

Assume that instead of the master system output $y$, we receive
the signal (\ref{27}).

In this case the adjustable observer with augmented error includes the
observer
\begin{equation}
\dot{\hat x}=A\hat x+\ph_0(y_r) +b\ph(y_r)\trn\hat
\theta+k(y_r-\hat y),~~~~\hat y= c\trn\hat x,
\label{28}
\end{equation}
scheme of augmentation
\begin{equation}
\bar  e=y_r - \hat y+H(p)\big[\bar \ph(t)\trn\hat \theta\big]-\bar
\omega(t)\trn\hat \theta,
\label{29}
\end{equation}
where $\bar \ph(t)=\ph(y_r(t))$, $\bar  \omega(t)=H(p)[\bar  \ph(t)]$,
and
robustified algorithm of adaptation
\begin{equation}
\dot{\hat \theta}=\gamma\bar \omega(t)\bar  e-\alpha(\hat \theta)\hat
\theta,
\label{30}
\end{equation}
\noindent where the function $\alpha(\hat \theta)$ obeys the
fol\-lo\-wing re\-lations
\begin{equation}
\alpha(\hat{\theta})=
\left\{\begin{array}{ll}
0,&|\hat{\theta}|<\theta^{*},\cr
\bigg(\dfrac{|\hat{\theta}|}{\theta^{*}}-1\bigg),&\theta^{*}\leq
|\hat{\theta}|\leq2\theta^{*},\cr
1,&|\hat{\theta}|>2\theta^{*}
\end{array}
\right.
\label{31}
\end{equation}
with any positive constant $\theta^{*}$.

\begin{theorem}\label{th:4}
The closed-loop system consisting of the master system
(\ref{2}), (\ref{27}), adjustable observer (\ref{28}),
scheme of augmentation (\ref{29}) and algorithm
of adaptation (\ref{30}), (\ref{31}) has the following properties:
\par{\it i}) for any initial conditions and any $\gamma\!>\!0$,
$\theta^*\!>\!0$ all the closed-loop signals are bounded and the
parameter error
$\tilde \theta$ converges
to the residual set
\begin{equation}
 D=\Big\{\tilde \theta:
|\tilde \theta|^2\le \max\big[(|\theta|+2\theta^*)^2,
\gamma\|\xi+\xi_e\|^2_\infty+|\theta|^2\big]\Big\},
\label{32}
\end{equation}
where the bounded variable $\xi_e$ obeys the equations:
\begin{equation}
\begin{aligned}
\dot\Xi_e&=F\Xi_e+\ph_0(y)-\ph_0(y_r)\\ 
&+b\Big(\ph(y)-\ph(y_r)
\Big)\trn\theta -k\xi,\\& \xi_e=c\trn\Xi_e;\label{nB}
\end{aligned}
\end{equation}

\par{\it ii}) if $\xi(t)\equiv 0$ and $\theta^*>|\theta|$ then, in
addition
to ({\it i}), regulation (\ref{u1}) is guaranteed.

\par{\it iii}) if, additionally, the vector-function $\ph(t)$
satisfies PE condition and the transfer function $H(p)$ is minimum
phase, then
convergence \eqref{u2} is achieved.
\end{theorem}
\begin{proof}
Differentiating estimation error $\eps=x-\hat x$ in view of equations
\eqref{2},
\eqref{27} and \eqref{28} after simple calculations we obtain:
\begin{equation}
\label{A1}
\dot\eps=F\eps+b\bar\ph(t)\trn\tilde\theta+\Delta(t),
\end{equation}
where $\Delta(t)=\ph_0(y)-\ph_0(y_r)+
b\big(\ph(y)-\ph(y_r)\big)\trn\theta
-k\xi$.

Then the augmented error $\bar e$ defined by equation \eqref{29} takes
the form
\begin{equation}
\label{A2}
\bar e=\bar\omega(t)\trn\tilde\theta+\xi_e+\xi,
\end{equation}
where the bounded variable $\xi_e$ obeys the equations (\ref{nB}).

Choose the Lyapunov function
$V(\tilde\theta)=\dfrac{1}{2}\tilde{\theta}\trn
\tilde\theta$. Its time derivative in view of \eqref{30} and
\eqref{A2} takes the form:
\begin{align*}
\dot V(\tilde\theta)&=\tilde{\theta}\trn\big(-\gamma\bar\omega\bar e
-\alpha\tilde\theta+\alpha \theta\big)\\
&=\tilde\theta\trn\big(-\gamma\bar\omega\bar{\omega}\trn\tilde{\theta}-
\gamma\bar\omega(\xi_e+\xi)-\alpha \tilde\theta+\alpha \theta\big) \\
&\le -\gamma\big|\bar{\omega}\trn\tilde\theta\big|^2+
-\sigma\big|\tilde\theta\big|^2+
\gamma\big|\bar{\omega}\trn\tilde\theta\big|\bigl\|\xi_e
+\xi\bigr\|_\infty+\alpha \big|\tilde\theta\big|\,
\big|\theta\big|\\
&\le -\dfrac{1}{2}\gamma\big|\bar{\omega}\trn\tilde\theta\big|^2-
\dfrac{1}{2}\alpha \big|\tilde\theta\big|^2+
\dfrac{1}{2}\gamma\left\|\xi_e+
\xi\right\|_\infty^2+\dfrac{1}{2}\alpha |\theta|^2\\
&\le -\dfrac{1}{2}\alpha \big|\tilde\theta\big|^2+
\dfrac{1}{2}\gamma\bigl\|\xi_e+\xi\bigr\|_\infty+\dfrac{1}{2}\alpha |\theta|^2.
\end{align*}
The latter inequality proves boundedness of all the closed-loop signals
and validity
of the estimate \eqref{32}.

If $\xi(t)\le 0$ and $\theta^*>|\theta|$, then the time derivative of
the Lyapunov function $V(\tilde\theta)=\dfrac{1}{2}\tilde{\theta}\trn \tilde\theta $ obeys the following expressions
\begin{align*}
\dot V(\tilde\theta)=-\gamma\big|\bar{\omega}\trn\tilde\theta\big|^2
+\gamma\sigma(\hat\theta)\tilde{\theta}\trn\hat\theta\le-
\gamma\big|\bar{\omega}\trn\tilde\theta\big|^2.
\end{align*}
The latter means validity of \eqref{u1}.
\end{proof}

To extend the above results to the master system (\ref{20}) with noisy
output (\ref{27}), we use the adjustable observer of the form
\begin{eqnarray}
\dot{\hat x}&=A(y_r)\hat x+\ph_0(y_r)+b\ph^T(y_r)\hat
\theta+k(y_r)(y_r-\hat y), \label{312}
\\ &\qquad \hat y=c\trn\hat x, \nonumber 
\end{eqnarray}
scheme of augmentation
\begin{eqnarray}
\bar e&=&y_r - \hat y+c^T\bar \eta \label{nC}\\
\dot{\bar  \eta}&=&\bar  G(t)\bar \eta - \bar \Omega\dot{\hat \theta},
\quad\bar \eta\in\mR^n,
\label{322}\\
\dot{\bar  \Omega} &=& \bar  G(t)\bar \Omega+b\bar \ph^T(t),
\quad\bar \Omega\in\mR^{n\times n}
\label{nD}
\end{eqnarray}
where $\bar  G(t)=A(y_r(t))-k(y_r(t))c\trn$, and
robustified algorithm of adaptation
\begin{equation}
\dot{\hat \theta}=\gamma\bar \omega(t)\bar  e-\alpha(\hat \theta)\hat
\theta,
\label{34}
\end{equation}
where $\bar \omega(t)=c\trn\bar  \Omega(t)$
and the function $\alpha(\hat \theta)$ obeys
relations (\ref{31}).

\begin{theorem}\label{th:5}
The closed-loop system consisting of the master system
(\ref{20}), (\ref{27}), adjustable observer (\ref{312}),
scheme of augmentation (\ref{nC})--(\ref{nD}) and algorithm
of adaptation (\ref{34}), (\ref{31}) has the following properties:
\par{\it i}) for any initial conditions and any $\gamma>0$,
$\theta^*>0$ all
the closed-loop signals are bounded and the parameter error $\tilde
\theta$
converges to the residual set
\begin{equation}
\begin{array}{l} 
 D=\Big\{\tilde \theta:|\tilde \theta|^2\le
\max\big[(|\theta|+2\theta^*)^2,
\gamma\|\xi+\xi_y\|^2_\infty+|\theta|^2\big]\Big\}
\label{332}
\end{array}
\end{equation}
where the bounded variable $\xi_y$ obeys the equations:
$$\begin{array}{l} 
\dot\Xi_y=\bar  G(t)\Xi_y+\ph_0(y)-\ph_0(y_r)\\
\quad+b\Big(\ph(y) -\ph(y_r)
\Big)\trn\theta
-k(y_r)\xi+\Big(A(y) -A(y_r)\Big)x,\\
\quad\xi_y=c\trn\Xi_y;
\end{array}$$

\par{\it ii}) if $\xi(t)\equiv 0$ and $\theta^*>|\theta|$ then, in
addition
to ({\it i}), regulation (\ref{u1}) is guaranteed.

\par{\it iii}) if, additionally, the vector-function $\ph(t)$
satisfies PE condition and the transfer function $H(p)$ is minimum
phase, then
convergence \eqref{u2} is achieved.
\end{theorem}
\begin{proof}
Differentiating estimation error $e=x-\hat x$ in view of equations
\eqref{20}, \eqref{27} and \eqref{312} after simple calculations we obtain:
\begin{equation}
\label{A3}
\dot \eps=\bar G(t)\eps+b\bar\ph(t)\trn\tilde\theta+\Delta(t),
\end{equation}
where $\Delta(t)=\big(A(y)-A(y_r)\big)x$ $+\ph_0(y)-\ph_0(y_r)$
$+b\big(\ph(y)-\ph(y_r)\big)\trn\theta$ $-k(y_r)\xi.$
Then differentiating the auxiliary error $\Xi_\delta=\eps+\bar\eta
-\bar\Omega\tilde\theta$ in view of equations \eqref{A3}, \eqref{332}
and \eqref{nD} we obtain
 $$
\begin{array}{l}
\dot\Xi_y=\bar  G(t)\Xi_y+\ph_0(y)-\ph_0(y_r)\\
\quad+b\Big(\ph(y) -\ph(y_r)
\Big)\trn\theta
-k(y_r)\xi+\Big(A(y)-A(y_r)\Big)x,\\
\quad\xi_y=c\trn\Xi_y;
\end{array}$$
Therefore, for the augmented error $\bar e$ defined by equation
\eqref{nC}
we can write $\bar e=\bar\omega\trn\tilde\theta+\xi_\delta+\xi,$ where
$\xi_\delta=c\trn\Xi_\delta$. Finally, using the same approach as in
the proof of Theorem \ref{th:4} we show validity of the all parts of 
Theorem \ref{th:5}.
\end{proof}

\section{Example: Signal Transmission via Adaptive Synchronization of
the Lorenz Systems}\label{Sec:4}

\subsection{\it Design of the aadptive observer}

{L}{et} us consider, for example, application of the proposed
method to adaptive synchronization of the {Lorenz systems}, exhibiting
chaotic behavior.

Let the master system be {\it Lorenz system}, described by the
following
equations
\cite{IEEE97,COS93}:
\begin{equation}
\left\{
\begin{array}{l}
\dot{x}_1=\sigma x_2-\sigma x_1,\\
\dot{x}_2=-x_2-x_1x_3+\theta x_1,\\
\dot{x}_3=-\beta x_3+x_1x_2.
\end{array}
\right.
\label{lor1}
\end{equation}
Constant parameters $\beta$, $\sigma$ are assumed to be known;
parameter
$\theta$ is varying depending on the information signal and its value
has to be
reconstructed by the observer. It is also assumed that the component
$x_1$ is
taken as a transmitted signal, i.e. $y\equiv x_1$.

Evidently, the system (\ref{lor1}) is a special case of (\ref{20}) with
the
following components:
\begin{equation}
\begin{aligned}
A(y)= \begin{bmatrix}-\sigma&\sigma&0\cr
0&-1&-y\cr
0&y&-\beta\end{bmatrix},~~ b=\begin{bmatrix}0\cr 1\cr 0\end{bmatrix},\\
\ph_0(y) \!=\!{\bf 0}_{3,1},~~ \ph(y) \!=\!y,~~ c\trn=[1,0,0].
\end{aligned}
\label{lor2}
\end{equation}
It is clear, that Assumption 4
is valid. To apply Theorem 3, accordingly with
Assumption 3, one has to find a vector-function $k(y)\inr^3$ so as the
system
$\dot x=\left(A(y)-k(y)c\trn\right)x,$ $y=c\trn x$
be asymptotically stable. Let us pick up $k(y)\equiv k=[0,-\sigma,
0]\trn$. Then
the matrix-function $G(y)= A(y)-k(y)c\trn $ is sum of a diagonal and a
skew-symmetric matrices:
\begin{equation}
G(y)= \begin{bmatrix}-\sigma&\sigma&0\cr
-\sigma &-1&-y\cr
0&y&-\beta\end{bmatrix}.
\label{lor3}
\end{equation}
It can be easily shown that this choice provides fulfillment of the
Assumption
3.
Indeed, let us consider the system $\dot x=G(c\trn x)x$, with the
matrix $G(y)$
given by (\ref{lor3}), and introduce the Lyapunov function $V(x)=0.5
x\trn x.$
Differentiating $V\big(x(t)\big)$ on $t$ one obtains:
\begin{align*}
\dot V(x)&=0.5 x\trn\big(
G(c\trn x)\trn+ G(c\trn x)\big)x\\
&=x\trn\begin{bmatrix}-\sigma&0&0\cr
0&-1&0\cr 0&0&-\beta\end{bmatrix}x \\
&=-\sigma x_1^2-x_2^2-\beta x_3^2.
\end{align*}
 Then the {\it exponential stability} of the system $\dot x=G(c\trn
x)x$ immediately follows. Therefore, Assumptions 3,4 are valid and
Theorem
\ref{th:3} can be applied. The adjustable observer for Lorenz-based
master
system (\ref{lor1}) has a form (\ref{22}), where the matrices are given
in
(\ref{lor2}). For estimation the unknown plant parameter $\theta$ (and,
thereby,
for recovering the message signal), the tuning algorithm
(\ref{23})--(\ref{26})
with $n=3$ has to be implemented in the observer. In the case of
significant
magnitude of the channel noise the robustified adaptation law
(\ref{30}),
(\ref{31}) can be used.

Some numerical examples of implementation of the proposed method for
adaptive synchronization of Lorenz systems are given below.

\subsection{\it Numerical example: square vaweform recovering via
adaptive synchronization of Lorenz systems}

Let us use the algorithm (\ref{22}), (\ref{23})--(\ref{26}) for
recovering parameter
$\theta$ of the master system (\ref{lor1}). It was assumed above that
$\theta$ is an unknown constant. In practice this is a varying information
signal, $\theta=\theta(t)$, and applicability of the proposed method depends on
the rate of tuning of the observer parameter $\hat\theta(t)$. This rate can be
found by means of the numerical examinations.

Let us rewrite the master/slave systems equations as follows:

\noindent{\bf Master system:}
\begin{equation}
\begin{array}{l}
\left\{
\begin{array}{l}
\dot{x}_1=\sigma x_2-\sigma x_1,\\
\dot{x}_2=-x_2-x_1x_3\\
\qquad\qquad+ r\big(1+\tet(t)\big) x_1,\\
\dot{x}_3=-\beta x_3+x_1x_2,
\end{array}
\right. \vspace{1mm}\\
\quad y(t)=x_1(t), \vspace{-2mm}
\end{array}
\label{lor1a}
\end{equation}
where $r$ is some known constant factor, $\tet(t)$ is
a varying parameter (in the case of communication via chaotic signal
modulation, $\tet(t)$ can be treated as an {\it information signal}).
It is fulfilled that $\theta= r(1+\tet)$.

\noindent {\bf Adjustable observer:}
\begin{equation}
\begin{array}{l}
\left\{
\begin{array}{l}
\dot{\hat x}_1=\sigma \hat x_2-\sigma \hat x_1,\\
\dot{\hat x}_2=-\hat {x}_2-y_r(t)\hat x_3+
\sigma e(t)\\
\qquad+ r\big(1+\hat\tet(t)\big) y_r(t),\\
\dot{\hat x}_3=-\beta x_3+y_r(t)\hat{x}_2,
\end{array}
\right. \vspace{1mm}\\
\quad e(t)=y_r(t)-\hat {x}_1(t), \vspace{-2mm}
\end{array}\label{lor2a}
\end{equation}
where $e(t)$ can be referred to as an {\it observation error}, $y_r(t)$ is
a measurable signal (in the case of communication systems, $y_r$ is
referred to as a {\it received signal}). For the noiseless case it is valid
that $y_r(t)\equiv y(t)$.

\noindent {\bf Augmented error filters:}
\begin{equation}
\begin{array}{l}
\left\{
\begin{array}{l}
\dot{\Omega}_1=\sigma \Omega _2-\sigma \Omega _1,\\
\dot{\Omega}_2=-\sigma \Omega _1-\Omega _2+y_r(t) \Omega _3\\
\dot{\Omega}_3=-\beta \Omega _3+y_r(t) \Omega _2,
\end{array}
\right. \vspace{1mm}\\
\quad \omega(t)=\Omega_1(t), \vspace{-2mm}
\end{array}
\label{lor3a}
\end{equation}

\begin{equation}
\left\{
\begin{array}{l}
\dot{\eta}_1=\sigma \eta_2-\sigma
\eta_1-\Omega_1(t)\dot{\hat\tet}(t),\\
\dot{\eta}_2=-\sigma \eta_1-\eta_2+y_r(t) \eta_3\\
\qquad\qquad -\Omega_2(t)\dot{\hat\tet}(t),\\
\dot{\eta}_3=-\beta \eta_3+y_r(t) \eta_2-\Omega_3(t)\dot{\hat\tet}(t).
\end{array}
\right.
\label{lor4a}
\end{equation}

\noindent {\bf Adaptation algorithm:}
\begin{equation}
\begin{array}{l}
\hat e(t)=e(t)+c\trn \eta(t),\\
\dot{\hat\tet}=\gamma\omega\hat e,\quad \hat\tet(0)=\hat\tet_0,
\end{array}
\label{lor5a}
\end{equation}
where parameter $\gamma>0$ is an {\it adaptation gain}.

The following numerical values of the master/slave systems parameters
were
chosen:
$$\sigma=10,~ \beta=8/3,~ r=97; ~~ \gamma=0.45~.$$

In our experiments, the ``square wave'' process $\tet(t)$ had been
taken.
In Fig. 1~{\it a} the time history of the measured signal
$y_r(t)$ is
shown. In the Fig. 1~{\it b} the time histories of
the original square waveform $\tet(t)$ and the waveform, recovered by
means of the adaptive observer  (\ref{lor2a})--(\ref{lor5a}) are shown. 

Figure~2 shows the corresponding time histories of the second
and third components of the state estimation errors $\eps_i(t)=x_i(t)$,
$i=1$,~$2$. Similar results for analogous information signal are shown
in Figs.~3, 4.  The simulation results demonstrate high adaptation rate of the 
proposed algorithm.
 It can be seen that the transient time for state estimates is the same as the 
one for parameter estimates.

%\begin{figure}[tbph]
%\begin{center}
%\unitlength=1mm
%\begin{picture}(120,112)
%\put(-15,112){\special{em:graph ythet.pcx}}
%\end{picture}
%\caption{{\it a}) Measured signal $y_r(t)$; {\it b})  square waveforms:
%$\tet$ -- original, $\hat\tet$ -- recovered.
%Algorithm \protect(\ref{lor2a})--(\protect\ref{lor5a}). }
%\end{center}
%\label{fig:1}
%\end{figure}

\begin{figure}[tbph]
\begin{center}
\includegraphics[width= 224pt,height= 168pt]{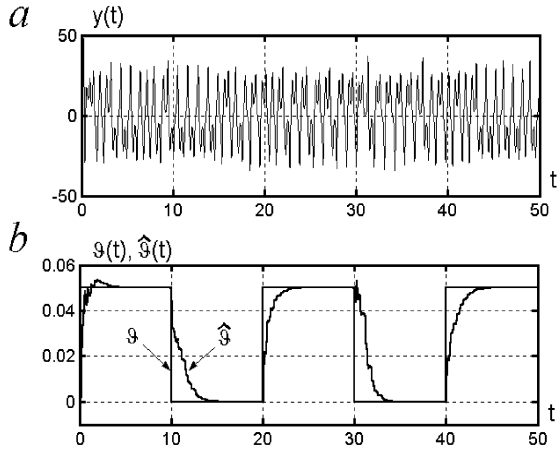}\\
\end{center}
\vspace{-0.3cm}
\caption{{\it a}) Measured signal $y_r(t)$; {\it b})  square waveforms:
$\tet$ -- original, $\hat\tet$ -- recovered.
Algorithm \protect(\ref{lor2a})--(\protect\ref{lor5a}). }
\label{fig:1}
\end{figure}
\begin{figure}[tbph]
\begin{center}
\includegraphics[width= 224pt,height= 168pt]{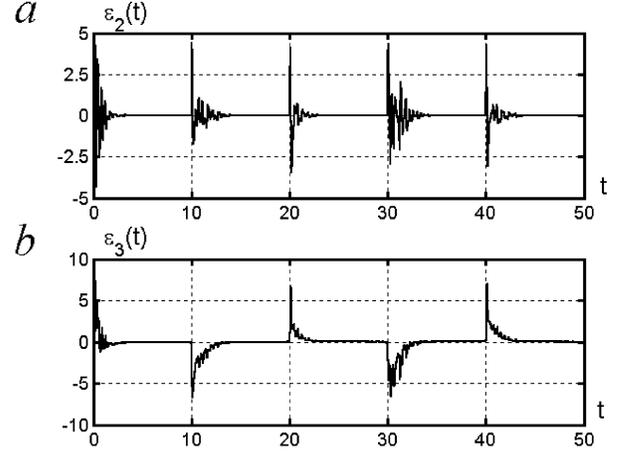}\\
\end{center}
\vspace{-0.3cm}
\caption{State estimation errors $\eps_2$, $\eps_3$ time histories.
Algorithm \protect(\ref{lor2a})--(\protect\ref{lor5a}).}
\label{fig:2}
\end{figure}

%\begin{figure}[tbph]
%\begin{center}
%\unitlength=1mm
%\begin{picture}(120,112)
%\put(-15,112){\special{em:graph eps23.pcx}}
%\end{picture}
%\caption{ State estimation errors $\eps_2$, $\eps_3$ time histories.
%Algorithm \protect(\ref{lor2a})--(\protect\ref{lor5a}).}
%\end{center}
%\label{fig:2}
%\end{figure}

%\begin{figure}[tbph]
%\begin{center}
%\unitlength=1mm
%\begin{picture}(120,112)
%\put(-15,112){\special{em:graph aythet.pcx}}
%\end{picture}
%\caption{{\it a}) Measured signal $y_r(t)$; {\it b}) analogous 
%information signals:
%$\tet$ -- original, $\hat\tet$ -- recovered.
%Algorithm \protect(\ref{lor2a})--(\protect\ref{lor5a}).}
%\end{center}
%\label{fig:3}
%\end{figure}

\begin{figure}[tbph]
\begin{center}
\includegraphics[width= 224pt,height= 168pt]{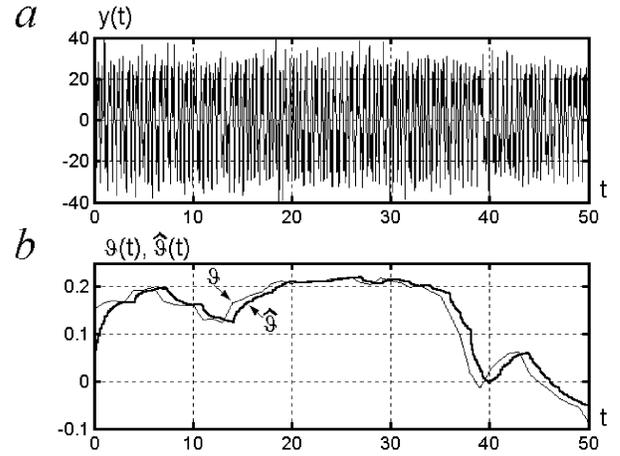}\\
\end{center}
\vspace{-0.3cm}
\caption{{\it a}) Measured signal $y_r(t)$; {\it b}) analogous 
information signals:
$\tet$ -- original, $\hat\tet$ -- recovered. 
Algorithm \protect(\ref{lor2a})--(\protect\ref{lor5a}).}
\label{fig:3}
\end{figure}
\begin{figure}[tbph]
\begin{center}
\includegraphics[width= 224pt,height= 168pt]{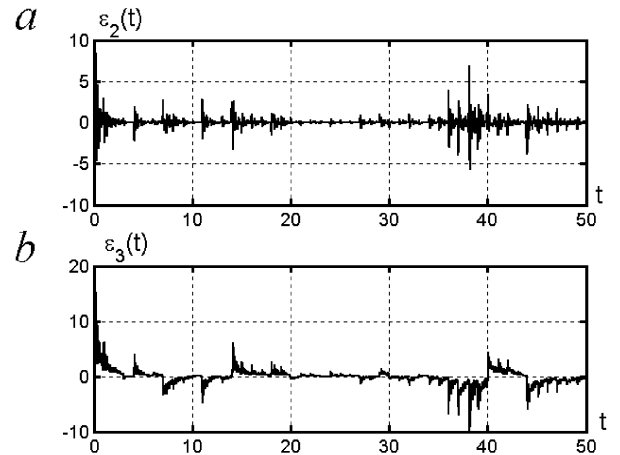}\\
\end{center}
\vspace{-0.3cm}
\caption{State estimation errors $\eps_2$, $\eps_3$ time histories.
Algorithm \protect(\ref{lor2a})--(\protect\ref{lor5a}), analogous 
information signal.} 
\label{fig:4}
\end{figure}

%\begin{figure}[tbph]
%\begin{center}
%\unitlength=1mm
%\begin{picture}(120,112)
%\put(-15,112){\special{em:graph aeps23.pcx}}
%\end{picture}
%\caption{State estimation errors $\eps_2$, $\eps_3$ time histories.
%Algorithm \protect(\ref{lor2a})--(\protect\ref{lor5a}), analogous 
%information signal.}
%\end{center}
%\label{fig:4}
%\end{figure}

\section{Conclusions}

{I}{n} this paper an unified approach for nonlinear adaptive
synchronization is proposed based on scheme of
\cite{F95,FM97,FNM00,AF00}
and new class of adaptation algorithms of \cite{FMN99,NikAut,NV01}.
It allows one to use chaotic signals generated by nonpassifiable
nonlinear
systems, particularly, by systems with relative degree greater than one
which
potentially increases security of communications.
Two versions of synchronization scheme are proposed based on augmented
error
adaptive observer and high-order tuners. Conditions of the estimate
convergence
are established for the noiseless case (Theorems
\ref{th:1}--\ref{th:3}).
Robustness of the scheme to the bounded measurement noise is
established
(Theorems \ref{th:4}--\ref{th:5}). The proposed approach allows one to
cope
with boun\-ded noise by means of robust modification of adaptation
algorithm.

Theoretical results are illustrated by simulation 
example by signal transmission based on adaptive synchronization of Lorenz
systems. That example demonstrates high of parameter identification
rate. The proposed algorithms may be applied to transmission of both binary 
(digital) and analogous signals in communication systems.

\nocite{*}
\bibliographystyle{IEEE}

\end{document}